\documentclass[11pt]{article}
\usepackage{amsfonts,amssymb, graphicx,a4,bm}
\usepackage{color}
\usepackage[english]{babel}
\usepackage{amsmath,amsthm,epsfig}
\usepackage[normalem]{ulem}

\newtheorem{theorem}{Theorem}[section]

\newtheorem{corollary}{Corollary}[section]
\newtheorem{definition}{Definition}[section]

\definecolor{GreenYellow}{cmyk}{0.15,0,0.69,0}
\definecolor{Yellow}{cmyk}{0,0,1,0}
\definecolor{Goldenrod}{cmyk}{0,0.10,0.84,0}
\definecolor{Dandelion}{cmyk}{0,0.29,0.84,0}
\definecolor{Apricot}{cmyk}{0,0.32,0.52,0}
\definecolor{Peach}{cmyk}{0,0.50,0.70,0}
\definecolor{Melon}{cmyk}{0,0.46,0.50,0}
\definecolor{YellowOrange}{cmyk}{0,0.42,1,0}
\definecolor{Orange}{cmyk}{0,0.61,0.87,0}
\definecolor{BurntOrange}{cmyk}{0,0.51,1,0}
\definecolor{Bittersweet}{cmyk}{0,0.75,1,0.24}
\definecolor{RedOrange}{cmyk}{0,0.77,0.87,0}
\definecolor{Mahogany}{cmyk}{0,0.85,0.87,0.35}
\definecolor{Maroon}{cmyk}{0,0.87,0.68,0.32}
\definecolor{BrickRed}{cmyk}{0,0.89,0.94,0.28}
\definecolor{Red}{cmyk}{0,1,1,0}
\definecolor{OrangeRed}{cmyk}{0,1,0.50,0}
\definecolor{RubineRed}{cmyk}{0,1,0.13,0}
\definecolor{WildStrawberry}{cmyk}{0,0.96,0.39,0}
\definecolor{Salmon}{cmyk}{0,0.53,0.38,0}
\definecolor{CarnationPink}{cmyk}{0,0.63,0,0}
\definecolor{Magenta}{cmyk}{0,1,0,0}
\definecolor{VioletRed}{cmyk}{0,0.81,0,0}
\definecolor{Rhodamine}{cmyk}{0,0.82,0,0}
\definecolor{Mulberry}{cmyk}{0.34,0.90,0,0.02}
\definecolor{RedViolet}{cmyk}{0.07,0.90,0,0.34}
\definecolor{Fuchsia}{cmyk}{0.47,0.91,0,0.08}
\definecolor{Lavender}{cmyk}{0,0.48,0,0}
\definecolor{Thistle}{cmyk}{0.12,0.59,0,0}
\definecolor{Orchid}{cmyk}{0.32,0.64,0,0}
\definecolor{DarkOrchid}{cmyk}{0.40,0.80,0.20,0}
\definecolor{Purple}{cmyk}{0.45,0.86,0,0}
\definecolor{Plum}{cmyk}{0.50,1,0,0}
\definecolor{Violet}{cmyk}{0.79,0.88,0,0}
\definecolor{RoyalPurple}{cmyk}{0.75,0.90,0,0}
\definecolor{BlueViolet}{cmyk}{0.86,0.91,0,0.04}
\definecolor{Periwinkle}{cmyk}{0.57,0.55,0,0}
\definecolor{CadetBlue}{cmyk}{0.62,0.57,0.23,0}
\definecolor{CornflowerBlue}{cmyk}{0.65,0.13,0,0}
\definecolor{MidnightBlue}{cmyk}{0.98,0.13,0,0.43}
\definecolor{NavyBlue}{cmyk}{0.94,0.54,0,0}
\definecolor{RoyalBlue}{cmyk}{1,0.50,0,0}
\definecolor{Blue}{cmyk}{1,1,0,0}
\definecolor{Cerulean}{cmyk}{0.94,0.11,0,0}
\definecolor{Cyan}{cmyk}{1,0,0,0}
\definecolor{ProcessBlue}{cmyk}{0.96,0,0,0}
\definecolor{SkyBlue}{cmyk}{0.62,0,0.12,0}
\definecolor{Turquoise}{cmyk}{0.85,0,0.20,0}
\definecolor{TealBlue}{cmyk}{0.86,0,0.34,0.02}
\definecolor{Aquamarine}{cmyk}{0.82,0,0.30,0}
\definecolor{BlueGreen}{cmyk}{0.85,0,0.33,0}
\definecolor{Emerald}{cmyk}{1,0,0.50,0}
\definecolor{JungleGreen}{cmyk}{0.99,0,0.52,0}
\definecolor{SeaGreen}{cmyk}{0.69,0,0.50,0}
\definecolor{Green}{cmyk}{1,0,1,0}
\definecolor{ForestGreen}{cmyk}{0.91,0,0.88,0.12}
\definecolor{PineGreen}{cmyk}{0.92,0,0.59,0.25}
\definecolor{LimeGreen}{cmyk}{0.50,0,1,0}
\definecolor{YellowGreen}{cmyk}{0.44,0,0.74,0}
\definecolor{SpringGreen}{cmyk}{0.26,0,0.76,0}
\definecolor{OliveGreen}{cmyk}{0.64,0,0.95,0.40}
\definecolor{RawSienna}{cmyk}{0,0.72,1,0.45}
\definecolor{Sepia}{cmyk}{0,0.83,1,0.70}
\definecolor{Brown}{cmyk}{0,0.81,1,0.60}
\definecolor{Tan}{cmyk}{0.14,0.42,0.56,0}
\definecolor{Gray}{cmyk}{0,0,0,0.50}
\definecolor{Black}{cmyk}{0,0,0,1}
\definecolor{White}{cmyk}{0,0,0,0}

\begin{document}
\def\blu{\color{Blue}}  
\def\mag{\color{Maroon}} 
\def\red{\color{Red}}  
\def\green{\color{ForestGreen}} 

\def\reff#1{(\protect\ref{#1})}
\let\a=\alpha \let\b=\beta \let\ch=\chi \let\d=\delta \let\e=\varepsilon
\let\f=\varphi \let\g=\gamma \let\h=\eta    \let\k=\kappa \let\l=\lambda
\let\m=\mu \let\n=\nu \let\o=\omega    \let\p=\pi \let\ph=\varphi
\let\r=\rho \let\s=\sigma \let\t=\tau \let\th=\vartheta
\let\y=\upsilon \let\x=\xi \let\z=\zeta
\let\D=\Delta \let\F=\Phi \let\G=\Gamma \let\L=\Lambda \let\Th=\Theta
\let\O=\Omega \let\P=\Pi \let\Ps=\Psi \let\Si=\Sigma \let\X=\Xi
\let\Y=\Upsilon
\let\e=\epsilon

\global\newcount\numsec\global\newcount\numfor
\gdef\profonditastruttura{\dp\strutbox}
\def\senondefinito#1{\expandafter\ifx\csname#1\endcsname\relax}
\def\SIA #1,#2,#3 {\senondefinito{#1#2}
\expandafter\xdef\csname #1#2\endcsname{#3} \else
\write16{???? il simbolo #2 e' gia' stato definito !!!!} \fi}
\def\etichetta(#1){(\veroparagrafo.\veraformula)
\SIA e,#1,(\veroparagrafo.\veraformula)
 \global\advance\numfor by 1
 \write16{ EQ \equ(#1) ha simbolo #1 }}
\def\etichettaa(#1){(A\veroparagrafo.\veraformula)
 \SIA e,#1,(A\veroparagrafo.\veraformula)
 \global\advance\numfor by 1\write16{ EQ \equ(#1) ha simbolo #1 }}
\def\BOZZA{\def\alato(##1){
 {\vtop to \profonditastruttura{\baselineskip
 \profonditastruttura\vss
 \rlap{\kern-\hsize\kern-1.2truecm{$\scriptstyle##1$}}}}}}
\def\alato(#1){}
\def\veroparagrafo{\number\numsec}\def\veraformula{\number\numfor}
\def\Eq(#1){\eqno{\etichetta(#1)\alato(#1)}}
\def\eq(#1){\etichetta(#1)\alato(#1)}
\def\Eqa(#1){\eqno{\etichettaa(#1)\alato(#1)}}
\def\eqa(#1){\etichettaa(#1)\alato(#1)}
\def\equ(#1){\senondefinito{e#1}$\clubsuit$#1\else\csname e#1\endcsname\fi}
\let\EQ=\Eq
\def\0{\emptyset}
\def\pp{{\bm p}}\def\pt{{\tilde{\bm p}}}

\def\\{\noindent}
\let\io=\infty
\def\VU{{\mathbb{V}}}
\def\EE{{\mathbb{E}}}
\def\GI{{\mathbb{G}}}
\def\TT{{\mathbb{T}}}
\def\C{\mathbb{C}}
\def\CC{{\mathcal C}}
\def\II{{\mathcal I}}
\def\LL{{\cal L}}
\def\RR{{\cal R}}
\def\SS{{\cal S}}
\def\NN{{\cal N}}
\def\HH{{\cal H}}
\def\GG{{\cal G}}
\def\PP{{\cal P}}
\def\AA{{\cal A}}
\def\BB{{\cal B}}
\def\FF{{\cal F}}
\def\v{\vskip.1cm}
\def\vv{\vskip.2cm}
\def\gt{{\tilde\g}}
\def\E{{\mathcal E} }
\def\I{{\rm I}}
\def\rfp{R^{*}}
\def\rd{R^{^{_{\rm D}}}}
\def\ffp{\varphi^{*}}
\def\ffpt{\widetilde\varphi^{*}}
\def\fd{\varphi^{^{_{\rm D}}}}
\def\fdt{\widetilde\varphi^{^{_{\rm D}}}}
\def\pfp{\Pi^{*}}
\def\pd{\Pi^{^{_{\rm D}}}}
\def\pbfp{\Pi^{*}}
\def\fbfp{{\bm\varphi}^{*}}
\def\fbd{{\bm\varphi}^{^{_{\rm D}}}}
\def\rfpt{{\widetilde R}^{*}}
\def\tende#1{\vtop{\ialign{##\crcr\rightarrowfill\crcr
              \noalign{\kern-1pt\nointerlineskip}
              \hskip3.pt${\scriptstyle #1}$\hskip3.pt\crcr}}}
\def\otto{{\kern-1.truept\leftarrow\kern-5.truept\to\kern-1.truept}}
\def\arm{{}}
\font\bigfnt=cmbx10 scaled\magstep1
\newcommand{\card}[1]{\left|#1\right|}
\newcommand{\und}[1]{\underline{#1}}
\def\1{\rlap{\mbox{\small\rm 1}}\kern.15em 1}
\def\ind#1{\1_{\{#1\}}}
\def\bydef{:=}
\def\defby{=:}
\def\buildd#1#2{\mathrel{\mathop{\kern 0pt#1}\limits_{#2}}}
\def\card#1{\left|#1\right|}
\def\proof{\noindent{\bf Proof. }}
\def\qed{ \square}
\def\trp{\mathbb{T}}
\def\trt{\mathcal{T}}
\def\Z{\mathbb{Z}}
\def\be{\begin{equation}}
\def\ee{\end{equation}}
\def\bea{\begin{eqnarray}}
\def\eea{\end{eqnarray}}

\title {A-priori Upper Bounds for the Set Covering Problem}
\author{Giovanni Felici$^{1}$, Sokol Ndreca$^{2}$,  Aldo Procacci$^3$ and Benedetto Scoppola$^4$ \\\\
\footnotesize{$^1$ Istituto di Analisi dei Sistemi ed Informatica, Consiglio Nazionale delle Ricerche,  00185 Roma}\\
\footnotesize{$^2$Dep. Estatistica-ICEx, UFMG, CP 702 Belo Horizonte - MG, 30161-970 Brazil}\\
\footnotesize{$^3$Dep. Matem\'atica-ICEx, UFMG, CP 702 Belo Horizonte - MG, 30161-970 Brazil}\\
\footnotesize{$^4$Dipartimento di Matematica - Universita Tor Vergata di Roma, 00133 Roma, Italy}\\
\tiny{emails: {giovanni.felici@iasi.cnr.it};~{sokol@est.ufmg.br};}\\
\tiny{aldo@mat.ufmg.br};~{scoppola@mat.uniroma2.it}}
\date{}
\maketitle
\numberwithin{equation}{section}
\begin{abstract}
\noindent
In this paper we present a new bound obtained with the probabilistic method for the solution of the Set Covering problem with unit costs. The bound is valid for problems of fixed dimension, thus extending previous similar asymptotic results, and it depends only on the number of rows of the coefficient matrix and the row densities. We also consider the particular case of matrices that are \textit{almost} block decomposable, and show how the bound may improve according to the particular decomposition adopted. Such final result may provide interesting indications for comparing different matrix decomposition strategies.
\end{abstract}
\section{Introduction}\label{intro}

\noindent
Given a finite ground set of objects $G$ and a finite collection $\mathcal{G}$ of its subset, a \textit{Set Cover} $\mathcal{C}$ is a subset of $\mathcal{G}$ such that each element of $G$ is contained in at least one of the subsets in $\mathcal{G}$.
\noindent
The \textit{Set Covering Problem} (SCP) consists in finding the set $\mathcal{C}$ of 
minimum cardinality. If positive weights are attached to each element of $\mathcal{G}$, the \textit{weighted} version of SCP  consists in finding a set $\mathcal{C}$ for which the sum of the weights of its element is minimum. Non weighted SCP may also be referred to as SCP with unit costs. {If $m$ and $n$ denote the cardinalities of $G$ and $\cal G$ respectively, the
SCP is usually reformulated as the problem of covering the rows of a $m\times n$ matrix $M$, whose
rows are associated with the $m$ elements of the ground set $G$, whose columns are associated with the $n$
subsets of $\cal G$, and whose entries are 1 if the element of the ground set associated with the row is contained  in the subset
associated with the columns, and 0 otherwise.}

\noindent
SCP is listed among the NP-complete problem class \cite{GA1979}, and is therefore considered to be a difficult problem to solve according to the fact that the solution time of any known algorithm cannot be bounded by a polynomial in the size the of the problem, unless of course, $P=NP$.

\noindent
Given their simplicity and generality, SCPs arise naturally in modeling many real-life problems, some interesting and large-sized examples of which can be found, among others, in
crew scheduling and allocation
(\cite{BURKE2014, CACCH2014, BOSCHE2014}),
data mining 
(\cite{LAD1988, LAD2005}), 
or (\cite{CHAO2007, LIA2009,VIJE2010}).

 Therefore, a large effort has been devoted by the research community to find efficient algorithms for its solution. Extensive discussion on algorithms to solve SCP can be found in several surveys, ranging from the 1975 Christofides  \cite{CH1975} to the more recent work of Fischetti et al.
 \cite{BO2000}.
An important stream of research on the issue is devoted to approximation results for SCP,
where solution algorithms are evaluated in their ability to find a solution whose distance from the optimum is guaranteed with a given probability.
Since SCP is NP-hard, numerous heuristic algorithms -- mostly of the greedy type -- have been developed for its solution, and the best approximation ratio available in polynomial time is
$H_d=\sum_{k=1}^{d}\frac{1}{k}$, i.e. $\Theta(\log d)$, where $d$ is the size
of the largest subset
(as from \cite{CHVA, fie1998, joh1974, LOV}), assuming $P\neq NP$ (recall that the approximation ratio of an algorithm is the ratio between the cost of the solution obtained by the algorithm and the cost of an optimal solution).

\noindent
More recently, a result by Levin \cite{levi2008} provides an approximation ratio of $H_d - \frac{196}{390}$, which improves the previous  results on the approximation ratio of the greedy algorithm. For the experimental results and comparison of the performance of many different  approximation algorithms of SCP, see e.g. \cite{gross1997, lan2007}.

\noindent
In 1984 Vercellis \cite{VE1984}, guided by previous results \cite{CHVA,LOV},
studied  SCP
problems with certain properties, namely  defined by a matrix $M$ with random i.i.d. Bernoulli entries.  In that paper
a class of randomized algorithms which find almost surely a solution whose approximation ratio tends asymptotically to 1 was exhibited, providing the asymptotic cardinality of  the optimal solution of SCP problems associated with random matrices in the above sense.

\noindent
In this paper we prove, via the so-called {\em probabilistic method in combinatorics} (see \cite{AS}), that the asymptotic SCP cardinality value for random matrices with fixed density $\d$ found in \cite{VE1984} is actually an upper bound valid for any matrix
with maximum row density $\d$ and any fixed dimension. Moreover, we show how such a-priori bound may be tailored for matrices with uneven row densities, and how the use of more refined approximations in the computations could result in a (sub-leading) improvement of its value.

\noindent
In addition, we define the class of the $(\nu,\mu)$-decomposable $0 - 1$ matrices and we show that when the  matrix of an SCP belongs to this class, an improved  bound can be obtained  according to the depth of the decomposition. Such fact indicates an interesting direction for applications, where SCP problems that are not perfectly decomposable may be treated with approximate decomposition algorithms guided by the evaluation of our bound.

\noindent
The paper is organized as follows: in Section \ref{notation} we provide the basic notation that will be used throughout the paper. Section \ref{bound} describes the main result and its extensions; Section \ref{decomp} discusses the refinement of the bound in the case of $(\nu,\mu)$-decomposable matrices. Some conclusions are drawn in Section \ref{conclu}.

\section{Notation and Previous Results}\label{notation}

A formal definition of SCP is given below:

\begin{definition} \label{scp}
Let $G=\{g_1, g_2, \dots, g_m\}$ be a ground set of $m$ elements, and let
$\mathcal{G}\subset 2^G$ be a collection of subset of $G$, $|\mathcal{G}|=n$,
where $\cup_{S\in \mathcal{G}}S=G$ and each $S$ has a positive cost $c_S$. $S_j\in {\cal G}$ covers $g_i\in G$ if $g_i\in S_j$.
$\mathcal{C} \subset \mathcal{G}$ is said to be a \textit{cover} of $G$ if $\cup_{S\in \mathcal{C}}S=G$.
A \textit{minimal cover} of $G$ is a cover for which $\sum_{S\in \mathcal{C}}c_S$ is minimum.
Given $G$, $\mathcal{G}$, and the associated costs $c_S, S \in \mathcal{G}$, the \textit{Set Covering Problem} SCP amounts to finding a minimal cover of $G$.
\end{definition}

\noindent
Without loss of generality, we assume that $G=\{1, 2, \dots, m\}$ and
$\mathcal{G}=\{S_j, j\in\{1, 2, \dots, n \} \}$.

\noindent
{ As said in the introduction}, we can describe SCP as the problem of covering the rows of a $m\times n$ matrix $M$, whose rows are associated with the $m$ elements of the ground set, whose columns are associated with the $n$
subsets of $\cal G$, and whose entries are 1 if the element of the ground set associated with the row is contained  in the subset
associated with the columns, and 0 otherwise. More formally

%

\begin{definition}\label{2.2}
Given a  $m\times n$ matrix $M=(m_{ij})$, where
\begin{equation*}
 m_{ij}=\begin{cases}
  1,   &\text{if column  $j$ with associated cost $c_j$ covers row $i$} \\
 0,   &\mbox{otherwise}\\
 \end{cases}
 \end{equation*}


\\then, SCP seeks the subset of columns that covers all rows, whose sum of costs is minimal.
When the cost $c_i=1$ for all $i$, the  problem is  called SCP with unit costs and the solution is given by the cover of minimal cardinality.
\end{definition}

\noindent
In this paper we deal with SCP with unit costs; in the following all SCPs are assumed to be of that type.

\vv
\\{\bf Remark}. Observe that if $M$ is a $0\!-\!1$ matrix describing the set covering problem of a given ground set  $G$ with
a given collection $\cal G$ of subsets of $G$, the any other matrix $M'$ obtained from $M$ by permutations of its rows and/or columns
describes the same set-covering problem as the original matrix $M$.

\vv

\noindent
Let us summarize briefly the results on SCP for random matrices. As remarked in the introduction, the first result,
obtained in \cite{VE1984}, is related to matrices in which $m_{ij}$ are  i.i.d. $0\!-\!1$ Bernoulli variable with probability $\d$. This model is known as the \textit{constant density model for SCP}.

\begin{theorem}[\bf Vercellis]
Let $C_m$ be the random variable that represent the optimal cost of random SCP. Suppose that the following two condition are satisfied:

\begin{itemize}
\item[C1]: $\qquad  \lim_{m\to \infty} \frac{n}{\log m}=\infty,$
\item [C2]: there exist $\alpha>0$ such that $n\le m^{\alpha}$.
\end{itemize}

\noindent
Then the sequence of random variables $C_m$ satisfies

\begin{equation}\label{verce}
 \lim_{m\to \infty} \frac{C_m}{\log m}=\Big[\log\frac{1}{1-\delta}\Big]^{-1} \quad a.s.
\end{equation}
\end{theorem}

\noindent

\\Note that, when $m$ and $n$ are asymptotically large, the optimal cost is given by
$$
C_m= \frac{\log m}{|\log(1-\d)|}
$$
with probability 1.

\noindent

For $i=\{1,2,\dots,m\}$, let  $\d_i$  be  the density  of $1$' of the row $i$ each row (simply called \textit{row density} from now on), i.e.,
$$
\d_i= {1\over n}\sum\limits_{j=1}^n m_{ij}
$$

A second model, introduced by Karp \cite{karp1976}, assumes that there is an equal number of ones in each row of the matrix $M$, that is, 
$$
\d_i= {1\over n}\sum\limits_{j=1}^n m_{ij}= \delta, \forall i
$$

\noindent
Note that in the Karp model the random variables $m_{ij}$ are not independent for $j\in\{1, 2, \dots, n\}$,
but they are indeed independent for $i\in\{1, 2, \dots, m\}$. These models have been studied  by Fontanari
\cite{font1996} using statistical mechanics techniques, which are useful  in the study of combinatorial optimization problems,
see e.g. Mezard, Parisi and Virasoro \cite{mezard1987}. The main result of Fontanari's work is that, for the Karp model,
the lower bound for the optimal cost is the same obtained by Vercellis  in the constant density model.

\section{An \textit{a-priori} Bound}\label{bound}

\noindent
Let $M$  be a given  $0\!-\!1$ matrix  with $m$ rows and $n$ columns.
Solving SCP for $M$ corresponds  to find a set  {$J\subset \{1,2,...,n\}$} of columns of $M$ of minimal cardinality $|J|$ such that for all $i\in\{1,2,...,m\}$
$$
\sum_{j\in J}m_{ij}>0
$$
\noindent
We are in particular interested in a possible a-priori estimate of the minimal cardinality $k=|J|$ of $J$ as a function of the  densities $\d_i$.  To get such upper  bound  we will use the so-called probabilistic method in combinatorics.

\noindent
The philosophy of the  probabilistic method is to prove the existence of combinatorial objects with certain desirable properties (e.g. a proper coloring of the edges of a graph) by showing that these objects have a positive probability to occur in some suitably defined probability space. In particular, the method works as follows. Suppose we are able to define
a probability space in which the occurrence of the combinatorial object with the desirable property -- the ``good event'' $A$ -- is ensured if a collection of ``bad events'' $\{B_1,\dots, B_m\}$ is such that none of them  occur. Namely we assume that we are able to define a probability space in which  the good event $A$ can be written as
$$
A = \bigcap_{i=1}^m \bar B_i
$$
where  $\bar B_i$ denote the {{probabilistic}} complement of $B_i$ (i.e. $\bar B_i$  is the event that $B_i$ does not occur). Suppose
then to be able to calculate (or to give an upper bound of) the probability $P(B_i)$ of occurrence for each of the bad events. Then,
the probability of the event $A$ is given by

$$
P(A)= P(\bigcap_{i=1}^m \bar B_i)= 1- P (\bigcup_{i=1}^m  B_i)
$$

\\Regardless the structure of dependencies of events $B_i$ we can write

\begin{equation}\label{unio}
P (\cup_{i=1}^m  B_i)\le \sum_{i=1} ^m P(B_i)
\end{equation}
\noindent

\noindent
{ Thus the good event $A$ occurs with positive probability if}

\begin{equation}\label{cond}
\sum_{i=1} ^m P(B_i) < 1
\end{equation}

\noindent
We note that inequality \ref{cond} is the well-known Local Lov\'asz Lemma condition (see, e.g., \cite{AS}) when, as it is in our case, each bad event depends on all the others.

\noindent
This philosophy can be applied to SCP for the fixed matrix $M$
in a quite straightforward way. Indeed, consider a probability space in which the elementary events are the uniformly random choices of a set $J$  with fixed cardinality $|J|=k$  of columns in the matrix $M$.  Define $m$ bad events $B_1,\dots, B_m$  with $B_i$ being the event that $\sum_{j in J}m_{ij}=0$. In other words, $B_i$ is the event that the $i$'th row is not covered by the columns in the set $J$. Then the good event $A$ is the event that ``{every row is covered by at least a column of the set $J$}'' and $A$ clearly occurs if none of the events $B_i$ occur. It is immediate to see that the probability $p_i=P(B_i)$ is such that

\begin{equation}\label{pai}
p_i\le(1-\d_i)^k
\end{equation}

\noindent
Indeed,
$$
p_i = \frac{{n-\d_in\choose k}}{{n\choose k}}=
(1-\d_i)^k\ \frac{1-\frac{1}{n(1-\d_i)}}{1-\frac{1}{n}}\; \frac{1-\frac{2}{n(1-\d_i)}}{1-\frac{2}{n}}\cdots
\frac{1-\frac{k-1}{n(1-\d_i)}}{1-\frac{k-1}{n}}\le(1-\d_i)^k
$$

\noindent

\noindent
{Now, using the condition (\ref{cond}), we have that a covering $J$ of cardinality $k$ exists if}

\begin{equation}\label{cond3}
\sum_{l\in\{1,2,...,m\}}(1-\d_l)^k < 1
\end{equation}

\noindent
Hence we have proved the following theorem:

\begin{theorem}\label{result}
Given the $m\times n$ matrix $M$ as defined above with density $\d_i$ for the $i$-th row, it always exists a covering $J$ of cardinality $k$ given by

\begin{equation}\label{cond4}
k=\min\{i\in \{1, \ldots, n\} | \sum_{l\in\{1,2,...,m\}}(1-\d_l)^i<1\}
\end{equation}
\end{theorem}

\noindent
Letting $\d=\max_i\{\d_i\}$ be the maximal row density of the matrix $M$, we get immediately the following corollary

\begin{corollary}\label{result2}
Given the $m\times n$ matrix $M$ defined above
if the density $\d_i$ for the $i$-th row does not exceed $\d$, then there  exists a covering $J$ of cardinality

\begin{equation}\label{condhom}
k > \frac{\log m}{|\log(1-\d)|}
\end{equation}
\end{corollary}

\vglue.5truecm
\noindent
{\bf Remark}.
One may ask how good are the bounds obtained by Theorem \ref{result} and Corollary \ref{result2}.
Recalling the result of \cite{VE1984}, and in particular  formula \ref{verce} in section \ref{notation},
we can observe that for matrices in which the only information available is the maximum row  density $\d$
 the bound (\ref{condhom}) is optimal in the sense that it is possible to exhibit an example of a matrix $M$ for which
  the optimal solution of the SCP has the cardinality given by the r.h.s. of (\ref{condhom}) asymptotically in $m, n$.
  Such matrix would, according to \cite{VE1984}, belong to the class of random matrices where entry is $1$ with probability
  $\d$ and $0$ otherwise. In other words, combining our result with that of \cite{VE1984},
  we can claim that the random matrices with constant density $\d$  have the worst possible
  optimal solution for the set covering problem (in the sense of the largest cardinality); i.e.,
  any other matrix $M$ with maximal (or even constant) row-density $\d$ has on optimal solution
  with cardinality less or equal than that of the random matrix with  density $\d$.

\vglue.5truecm

\subsection{Further refinements}

The above bounds (\ref{cond4}) and (\ref{condhom}) can be improved when additional information on
the structure of the matrix $M$, beside the row densities $\d_i$, is available.
This yields anyway into sub-leading corrections to the asymptotic bounds (\ref{cond4}) and (\ref{condhom}). The idea is the following. Starting from equation (\ref{unio}), we can give a better bound of the quantity $P(\cup_{i=1}^m B_i)$ { using the Bonferroni inequality}.
Indeed, instead of the trivial inequality (\ref{unio}), we can write
\begin{equation}\label{ext}
P(\cup_{i=1}^m B_i)\le\sum_{i=1}^{m}P(B_i)-\sum_{1\le i<j\le m}P(B_i\cap B_j)
+\sum_{1\le i<j<k\le m}P(B_i\cap B_j\cap B_k)
\end{equation}

\\The two probabilities $P(B_i\cap B_j)$ and $P(B_i\cap B_j\cap B_k)$ can be evaluated as follows. Let 
$\G_{ij}=\{l\in\{1,...,n\}|m_{il}=m_{jl}=1\}$. In other words, the set  $\G_{ij}$ represents the overlap between the two rows $i$ and $j$. We now define $\g_{ij} = \frac{|\G_{ij}|}{n}$.
Moreover, let $\G_{ijk}=\{l\in\{1,...,n\}| m_{il}=m_{jl}=m_{kl}=1\}$ and define
$\g_{ijk} = \frac{|\G_{ijk}|}{n}$.
Then, we can write

\begin{equation}\label{over2}
P(B_i\cap B_j)=\frac{{n(1-\d_i-\d_j+\g_{ij})\choose k}}{{n\choose k}}
\end{equation}

\begin{equation}\label{over3}
P(B_i\cap B_j\cap B_k)=\frac{{n(1-\d_i-\d_j-\d_k+\g_{ij}+
\g_{ik}+\g_{jk}-\g_{ijk})\choose k}}{{n\choose k}}
\end{equation}

\noindent
Then we can identify our upper bound on SCP finding the smallest $k$ such that
$$
\sum_{i=1}^{m}{n(1-\d_i)\choose k}-\sum_{1\le i<j\le m}{n(1-\d_i-\d_j+\g_{ij})\choose k}+
$$
\begin{equation}\label{condover}
+\sum_{1\le i<j<k\le m}{n(1-\d_i-\d_j-\d_k+\g_{ij}+
\g_{ik}+\g_{jk}-\g_{ijk})\choose k}<{n\choose k}
\end{equation}

\noindent
The above condition is easy to be checked numerically. Clearly there are choices of $\g_{ij}$ and $\g_{ijk}$ for which the condition above gives an estimate for $k$ which is better than (\ref{condhom}). To give a flavour, let us consider a random matrix with constant density $\d$. We have then $\g_{ij}=\d^2$ for all $1\le i<j\le m$ and $\g_{ijk}=\d^3$ for all $1\le i<j<k\le m$. Condition (\ref{condover}) becomes, neglecting $o(m)$ terms
$$
y-\frac{y^2}{2}+\frac{y^3}{6}<1
$$
where $y=m(1-\d)^k$. This gives, instead of (\ref{condhom}), the condition

$$
k>\frac{\log m-\log(1.56)}{|\log(1-\d)|}
$$

\noindent
In the (easy) latter case of the random matrix, the constant can be further improved by considering the intersections up to 5, 7, 9... sets $B_i$. This gives, in place of 1.56, larger and larger constants, and when we arrive up to the intersections of $m$ sets, for $n$ exponentially large in $m$, we  obtain $k > 0$. This is not surprising, since if the random matrix has a number of column exponentially large, then, almost surely, we  have a column $j$ with $m_{ij}=1\quad\forall i$.

\section{Improvements for partially decomposable problems}\label{decomp}

\noindent
The remark after Definition \ref{2.2} suggests that it may be convenient to permute rows
and columns of a $0\!-\!1$ matrix $M$ and to decompose it in blocks, in order to try to improve the bound
(\ref{condhom}). Indeed, the decomposition of a matrix is an alternative representation that allows a particular
structure to emerge. Such structures attract the interest of researchers as they may facilitate the solution of
certain mathematical problem where the decomposed matrix plays a role. 
It is well established (see e.g. \cite{BORNDORFER1998})
that decompositions allow to confine and control particularly ``difficult'' substructures
of the matrix, and moreover allow to parallelize solution algorithms over the substructures -- typically, \textit{blocks} --
identified by the decomposition. A comprehensive analysis of the different types of matrix decomposition and the related
algorithms is beyond the scope of this paper. Recently, similar problems are discussed in \cite{KAHOU2008} and \cite{BERTOLA2014};
the relevant interactions beween the block decomposition of binary matrix and several data mining problems are highlighted in \cite{TAO2005, ZHANG2005}.

\noindent
Here we define a class of decomposable $0\!-\!1$  matrices with maximum row density $\d$, that we assume to possess an interesting
structure, and exploit the extension of the bound (\ref{condhom}) for SCP whose associated matrix $M$ belongs to this class.
The general idea of this class is that the matrix can be  decomposed into four  block matrices such that the maximum row density
of the two block matrices in the main diagonal is larger that $\d$, while the maximum row density of the remaining off-diagonal
two matrices is smaller that $\d$. Such definition is indeed similar to the {\it bordered block diagonal form} treated in \cite{BORNDORFER1998},
where its interest for optimization problems is discussed and several decomposition algorithms are referred.

\begin{definition}\label{numudecomp}Let $M$ be a $0\!-\!1$  matrix of dimension $m \times n$ with maximum row density $\d$.
Let $\n,\m\in (-1,1)$. Then $M$ is \mbox{{\bf $(\nu,\mu)$-decomposable}} if, after a permutation of its rows and columns, becomes a $0-1$ matrix  $M'$  formed by  the 4 submatrices $M_{11}, M_{12}, M_{21}, M_{22}$ such that:

$$ M' =
\begin{pmatrix}
M_{11} & M_{12}\\
M_{21} & M_{22}
\end{pmatrix}
$$

and
\begin{itemize}{
\item $M_{11}$ has $\frac{m}{2}(1+\mu)$ rows,  $\frac{n}{2}(1+\nu)$ columns, and  maximum row density $\d_1>\d$,
\item $M_{12}$ has $\frac{m}{2}(1+\mu)$ rows,   $\frac{n}{2}(1-\nu)$ columns, and  maximum row density $\d_2<\d$,
\item $M_{21}$ has $\frac{m}{2}(1-\mu)$ rows, $\frac{n}{2}(1+\nu)$ columns, and  maximum row density
$\d_3<\d$,
\item $M_{22}$ has $\frac{m}{2}(1-\mu)$ rows,  $\frac{n}{2}(1-\nu)$ columns, and  maximum row density $\d_4>\d$.}
\end{itemize}
\end{definition}

\\We want to exploit the case of $(\nu,\mu)$-decomposable matrices with the following random experiment. Choose uniformly at random $k_1$ columns of the matrix $M$ in the first $\frac{n}{2}(1+\nu)$ columns, and $k_2$ columns in the following $\frac{n}{2}(1-\nu)$ columns.

%
%
%
%
%
Call $B_i$ the (bad) event to have that the $k_1+k_2$ columns do not cover the row $i$, with $i=1,2,...,\frac{m}{2}(1+\mu)$, and call $\widetilde B_{i}$ the (bad) event to have that the $k_1+k_2$ columns do not cover the row $i$, with $i=\frac{m}{2}(1+\mu)+1,...,m$. The probabilities of such bad events are, reasoning as in section \ref{bound}, bounded by:
\begin{equation}\label{probab}
P(B_i)\le(1-\d_1)^{k_1}(1-\d_2)^{k_2}\atop
P(\widetilde B_i)\le(1-\d_3)^{k_1}(1-\d_4)^{k_2}
\end{equation}

\\Calling $B=\cup_{i=1}^{\frac{m}{2}(1+\mu)} B_i$ and $\widetilde B=\cup_{i=\frac{m}{2}(1+\mu)+1}^m \widetilde B_i$,
we are looking for an a-priori estimate of the probability $P(B\cup \widetilde B)$ of the event $B\cup \widetilde B$ of the form

\begin{equation}\label{koli}
P(B\cup \widetilde B)=P(B)+P(\widetilde B)-P(B\cap \widetilde B)<1
\end{equation}

{
\\because this relation would imply, as before, the fact that the complementary event 
$A= \overline{B\cup \widetilde B}$, which is the (good) event to have all the rows covered by the  $k_1+k_2$ columns chosen by our random experiment, would have a probability strictly positive, and hence it would exist.}

\noindent
The trouble with (\ref{koli}) is the fact that it is not easy, in general, to give a non trivial lower bound of the quantity $P(B\cap \widetilde B)$.
Here we will use the trivial bound $P(B\cap \widetilde B)\ge 0$, and we will get rid of such term from the inequality.
There are, however, some specific cases in which it is easy to estimate that intersection: for instance, if $\d_2$ and $\d_3$ are zero
it is easy to see that $B$ and $\tilde B$ are independent, and then $P(B\cap \widetilde B)=P(B) \times P(\widetilde B)$. We will recall this later;
for the time being let us write our condition in terms of the following inequality:

\begin{equation}\label{nokoli}
P(B\cup \widetilde B)\le P(B)+P(\widetilde B)<1
\end{equation}
{
\\Plugging  (\ref{probab}) in (\ref{nokoli}), we get that if  $k_1$ and $k_2$ are such that
\begin{equation}\label{conk12}
\frac{m}{2}(1+\mu)(1-\d_1)^{k_1}(1-\d_2)^{k_2}+\frac{m}{2}(1-\mu)(1-\d_3)^{k_1}(1-\d_4)^{k_2}<1
\end{equation}

\\then  the good event $A= (B\cup \widetilde B)^ç$ to have all the rows covered by the  $k_1+k_2$ columns
has a positive probability to occur. The condition  (\ref{conk12}) can be separated in the two independent bounds:

\begin{equation}\label{sys1}
\frac{m}{2}(1+\mu)(1-\d_1)^{k_1}(1-\d_2)^{k_2}<\a
\atop
\frac{m}{2}(1-\mu)(1-\d_3)^{k_1}(1-\d_4)^{k_2}<1-\a
\end{equation}

\\with $\a\in (0,1)$ to be determined in order to optimize globally our bound. The two conditions in (\ref{sys1}) can be rewritten as

\begin{equation}\label{sys2}
k_1|\log(1-\d_1)|+k_2|\log(1-\d_2)|>c_1(\a)
\atop
k_1|\log(1-\d_3)|+k_2|\log(1-\d_4)|>c_2(\a)
\end{equation}
with
\begin{equation}\label{c1c2}
c_1(\a)=|\log\a|+\log \left[{m\over 2}(1+\m)\right]
\atop
c_2(\a)=|\log(1-\a)|+\log \left[{m\over 2}(1-\m)\right]
\end{equation}
Let us consider  the following linear system in $k_1$ and $k_2$
\begin{equation}\label{sys3}
\left\{k_1|\log(1-\d_1)|+k_2|\log(1-\d_2)|=c_1(\a)
\atop
k_1|\log(1-\d_3)|+k_2|\log(1-\d_4)|=c_2(\a)
\right.
\end{equation}
}

\\The solution of such a system is

\begin{equation}\label{k12}
k_1=\frac{1}{\D}\left( c_1(\a)|\log(1-\d_4)|-c_2(\a)|\log(1-\d_2)|\right)
\atop
k_2=\frac{1}{\D}\left( c_2(\a)|\log(1-\d_1)|-c_1(\a)|\log(1-\d_3)|\right)
\end{equation}
with $\D=\log(1-\d_1)\log(1-\d_4)-\log(1-\d_2)\log(1-\d_3)$.

\noindent
Now let us   find the $\a$ that minimizes the value of

\begin{equation}\label{sumk}
k_1+k_2=\frac{1}{\D}\left( c_1(\a)[|\log(1-\d_4)|-|\log(1-\d_3)|]+c_2(\a)[|\log(1-\d_1)|-|\log(1-\d_2)]|\right)
\end{equation}

\\{
It is easy to see that the unique minimum of $k_1+k_2$, as  $\a$  varies in $(0,1)$ is attained for $\a=\bar \a$, with $\bar\a$ given by
\begin{equation}\label{alp}
\bar\a=\frac{|\log(1-\d_4)|-|\log(1-\d_3)|}{|\log(1-\d_4)|-|\log(1-\d_3)|+|\log(1-\d_1)|-|\log(1-\d_2)}|
\end{equation}
Recalling now the explicit expression of $c_1(\a)$ and $c_2(\a)$ given in (\ref{c1c2}), we can put the value of $\a=\bar\a$  given by (\ref{alp}) into $c_1(\a),c_2(\a)$ appearing in (\ref{sumk}) and get an a priori upper bound for the cardinality $k=k_1+k_2$ of the optimal  solution in the case  of the decomposable matrix.
}

\\{Hence we have proved the following theorem.
\begin{theorem}\label{resultato2}
Let  $M$ be a $(\nu, \mu)$-decomposable matrix as given in definition \ref{numudecomp}.
Then there exists a covering $J$ of cardinality $k_1+k_2$ given by \eqref{sumk} with $\a$  given by \eqref{alp}.
\end{theorem}
}


\noindent
In general it is not simple to compare analytically the bound  (\ref{sumk}) (putting of course  $\a=\bar \a$ given by  (\ref{alp}))  with the bound (\ref{condhom}), but one can check numerically  that this estimate tends to  improve the previous
general estimate. In any case, it is thus well established that the relations among $\delta_1, \delta_2, \delta_3, \delta_4$ may play a role in the design of a solution algorithm.
We illustrate this fact by  considering   two  examples   in which the expression \eqref{sumk} simplifies drastically and yet are representative of possible situations.
For both these example  we  get an explicit improvement of the bound (\ref{condhom}).

\vv
\\{\bf Example 1}. Suppose that a proper  permutation of the rows and the columns of $M$ results in a perfect block decomposition, where, w.l.o.g.,
$$
\delta_2=\delta_3 = 0,~~~~
\delta_1 = {2\d\over 1+\n}, ~~~ \delta_4 = {2\d\over 1-\n} ~~~~~~~|\nu|<1-2\d
$$
Then the optimal solution of the SCP defined  by matrix $M$ can be obtained by the union of the solutions obtained on $M_{11}$ and $M_{22}$. Indeed in this case (i.e. $\d_2=\d_3=0$) the events $B$ and $\widetilde B$ are clearly independent. Therefore we can write the condition (\ref{koli}) in terms of

\begin{equation}\label{koli2}
P(B\cup \widetilde B)=P(B)+P(\widetilde B)-P(B)P(\widetilde B)<1
\end{equation}

\\and this is equivalent to impose separately $P(B)<1$ and $P(\widetilde B)<1$. It follows that in formulas \eqref{sys2}-\eqref{sumk} the factors  $c_1(\a)$ and $c_2(\a)$ can be replaced  by

\begin{equation}
c_1\doteq c_1(1)=\log\left[{m\over 2}(1+\mu)\right]
\atop
c_2\doteq c_2(1)=\log\left[{m\over 2}(1-\mu)\right]
\end{equation}

\\and  $k_1+k_2$ is such that

\begin{equation}\label{delta0}
k_1+k_2 = {\log\left[{m\over 2}(1+\mu)\right]\over |\log(1-\d_1)|}+  {\log\left[{m\over 2}(1-\mu)\right]\over |\log(1-\d_4)|}
\end{equation}
I.e., the resulting bound is exactly the general bound given in the section \ref{bound} applied separately to the two factorized problems in the blocks 1 and 4.

\\Let us now show that
$$
{\log\left[{m\over 2}(1+\mu)\right]\over |\log(1-\d_1)|}+  {\log\left[{m\over 2}(1-\mu)\right]\over |\log(1-\d_4)|} <  {\log m \over |\log(1-\d)|}
$$
which is equivalent to show that
$$
\log(m) \times
\left[{1\over{|\log(1-\d)|}} - {1\over{|\log(1-\d_1)|}} - {1\over{|\log(1-\d_4)|}}\right]
-
{{\log{{1+\mu} \over 2}} \over{|\log(1-\d_1)|}}
-
{{\log{{1-\mu} \over 2}} \over{|\log(1-\d_4)|}}> 0
$$
The last two terms do not depend on $m$ and are always non negative (recall that $|\m|<1$ and thus  $\log{{1\pm\mu} \over 2}<0$.
Therefore  it is enough to prove the
the following  inequality:

$$
\left[{1\over{|\log(1-\d)|}} - {1\over{|\log(1-\d_1)|}} - {1\over{|\log(1-\d_4)|}}\right]> 0
$$
i.e.
recalling  that $\d_1 = {{2\d} \over{1+\nu}}$ and $\d_4 = {{2\d} \over{1-\nu}}$, it is enough to prove
\be\label{ineq}
{1\over{|\log(1-\d)|}} >
{1\over{|\log(1-{{2\d}\over{1+\nu}})|}} +
{1\over{|\log(1-{{2\d}\over{1-\nu}})|}}
\ee

\\To show that the above inequality is always satisfied, first recall that we always have $|\nu| < 1- 2\d$.
One can now study the expression on the  r.h.s. of \ref{ineq}, as a function of $\nu\in  (-1+2\d,1-2\d)$.
Let
$$
f(\nu)= {1\over{|\log(1-{{2\d}\over{1+\nu}})|}} +
{1\over{|\log(1-{{2\d}\over{1-\nu}})|}}
$$
It can be checked that $f(\nu)$ is  concave  in the interval $\nu\in  (-1+2\d,1-2\d)$ and attains the maximum at $\nu=0$ where reaches the value
$$
f(0)= {2\over{|\log(1-{{2\d}})|}}
$$
and since
$$
{1\over{|\log(1-\d)|}} > {2\over{|\log(1-{{2\d}})|}}
$$
the inequality  above is true for all  $\nu\in  (-1+2\d,1-2\d)$.

\vv
\\{\bf Example 2}.  Suppose there exists a permutation of rows and  columns of the $m\times n$ matrix $M$ that identifies $4$ block matrices of size ${m\over 2}\times {n\over 2}$, such that the two blocks on the main diagonal absorb a large portion of the matrix density, at the expenses of the blocks in the other diagonal. Namely suppose that after such a permutation we have
$$
\d_1=\d_4=2\d-\e~~~~~~~~\d_2=\d_3=\epsilon~~~~~~~~\mu=\nu=0~~~~~~\e<\d
$$
In this case by \eqref{alp}  we get that $\bar\a={1\over 2}$ and thus, plugging this value in formulas \eqref{c1c2} we get
$$
c_1({1/2})= c_2({1/2})= \log m
$$
Therefore, considering that in the present case $\D= [\log(1-2\d+\e)]^2-[\log(1-\e)]^2$ equation \eqref{sumk} becomes
$$
k_1 + k_2 = {{2\log m} \over {|\log(1-2\d+\e)| + |\log(1-\e)|}}
$$
Let us prove that, for any $\e\in [0, \d)$ we have
$$
 {{2\log m} \over {|\log(1-2\d+\e)| + |\log(1-\e)|}}< \frac{\log m}{|\log(1-\d)|}
$$
which is equivalent to  the inequality

$$
|\log(1-2\d+\e)| + |\log(1-\e)|  > 2 |\log(1-\d)|
$$
i.e.
$$
\log\left[(1-2\d+\e)\times(1-\e)\right]  <\log\left[(1-\d)^2\right]
$$
Inequality above, by the monotonicity  of the logarithm,  is true if and only if
$$
(1-2\d+\e)\times(1-\e)  <  (1-\d)^2
$$
i.e. if
$$  2\d\e -\e^2   < \d^2$$
which is always true for all $\e\in [0,\d)$.

\section{Conclusions}\label{conclu}

The results of this paper are related with the existence of an easy to compute \textit{a-priori} upper bound for the Set Covering problem with unit cost. The bound is obtained by the application of the probabilistic method in combinatorics and extends  to a deterministic setting previous asymptotic results. We show several variants of the bound that can be computed by a simple binary search, and analyze some extensions. As a side results, we consider the specialization of this bound when the ${0-1}$ matrix that describes the SCP can be almost decomposed into a block diagonal matrix. In the latter case we show how the bound is related with the parameters that define the decomposition and show that, under certain conditions, the decomposition always improves the bound.

{
Although the results presented are mainly related with theoretical properties of the solution of a specific integer programming problem, we believe that they provide an interesting insight for practical application, given the extremely general and simple nature of the bound; moreover, the results of Section \ref{decomp} suggest that even non-perfect decompositions may be useful to improve solution methods for the hard combinatorial problems considered in this paper. Such considerations demand further investigations and computational tests that will be addressed in future research.
}

\section*{Acknowledgments}
This work  has been supported by Conselho Nacional de
Desenvolvimento Cient\'{i}fico e Tecnol\'ogico (CNPq), Consiglio Nazionale delle
Ricerche (CNR), and FAPEMIG
(Funda{c}\~ao de Amparo \`a Pesquisa do Estado de Minas
Gerais)-Programa Pesquisador Mineiro.


\begin{thebibliography}{99}

\small{
\bibitem{AS} N. Alon; J. Spencer (2008): {\it The Probabilistic Method.
Third  Edition}. New York, Wiley-Interscience.

\bibitem{BERTOLA2014} E. Bertolazzi; A. Rimoldi (2014): {\it  Fast matrix decomposition in $F^2$}, Journal of Computational and Applied
Mathematics, {\bf 260}, 519--532.

\bibitem{BORNDORFER1998} R. Borndorfer; C. E. Ferreira; A. Martin (1998): {\it Decomposing Matrices into Blocks},
Siam Journal of Optimization, {\bf 9}, n. 1, 236--269.

\bibitem{LAD2005} E. Boros; P.L. Hammer; T. Ibaraki (2005):
{\it Logical Analysis of Data}, 
In: Encyclopedia of Data Warehousing and Mining, (J. Wang, ed.) Idea Group Reference, 
689--692.

\bibitem{BOSCHE2014} M. Boschetti; V. Maniezzo (2014):
{\it A set covering based metheuristic for a real-world city logistics problem},
International Transactions in Operational Research, doi: 10.1111/itor.12110.

\bibitem{BURKE2014} E.K. Burke; T. Curtois (2014):
{\it New approaches to nurse rostering benchmark instances},
European Journal of Operational Research {\bf 237}, 71–-81.

\bibitem{CACCH2014} V. Cacchiani; V.C. Hemmelmayr; F. Tricoire (2014): 
{\it A set-covering based heuristic algorithm for the periodic vehicle routing problem},
Discrete Applied Mathematics, {\bf 163}, 53–-64.

\bibitem{BO2000} A. Caprara; P. Toth; M. Fischetti (2000): {\it Algorithms for the Set Covering Problem}, Annals of Operations Research, {\bf 98}, 353–-371.

\bibitem{CHAO2007} W.A. Chaovalitwongse; T.Y. Berger-Wolf; B. Dasgupta; M.V. Ashley (2007):
{\it Set covering approach for reconstruction of sibling relationships},
Optimization Methods and Software, {\bf 22}, 11–-24.

\bibitem{LIA2009} L. Chen; J. Crampton (2009): 
{\it Set Covering Problems in Role-Based Access Control},
Lecture Notes in Computer Science {\bf 5789}, 689--704.

\bibitem{CH1975} N. Christofides; S. Korman (1975): {\it  A Computational Survey of Methods for the Set Covering Problem}, Management Science, {\bf 21}, 591--599.

\bibitem{CHVA} V. Chvatal (1979): {\it A greedy heuristic for the set-covering problem}. Math. Oper. Res. {\bf 4}, no. 3, pp. 233--235.

\bibitem{LAD1988} Y. Crama; P.L. Hammer; T. Ibaraki (1988):
{\it Cause-e?ect relationships and partially de?ned Boolean functions},
Annals of Operational Research, {\bf 16}, 299–-325.

\bibitem{fie1998} U. A. Fiege (1998): {\it Threshold of ln n for approximating set cover}, Journal of the ACM, {\bf 45} (4),
634-652.

\bibitem{font1996}J. F. Fontanari (1996): {\it  A statistical mechanics analysis of the set covering problem}, J. Phys.
A: Math. Gen., {\bf 9}, 473--483.


\bibitem{GA1979}  M. R. Garey; D. S. Johnson (1979): {\it Computers and Intractability: A Guide to the Theory of NP-Completeness}, Freeman and Co.

\bibitem{gim1967} J. F. Gimpel (1967): {\it A Stochastic Approach to the Solution of Large Covering Problems},
IEEE Switching and Automata Theory, 76--83.

\bibitem{GOLD} O. Goldschmidt; D. S. Hochbaum; G. Yu (1993):    {\it A modified greedy heuristic for the set covering problem with improved worst case bound},
Information Processing Letters archive, 48-6, Dec. 20, 1993,
305--310.

\bibitem{gross1997} T. Grossman; A. Wool (1997):  {\it  Computational experience with approximation algorithms
for the set covering problem}, European Journal of Operational Research, {\bf 101}, 81--92.

\bibitem{joh1974} D. S. Johnson (1974):  {\it Approximation algorithms for combinatorial problems}, J. Comput.
System Sci., {\it 9}, 256--278.

\bibitem{KAHOU2008} G. A. A.  Kahou; L. Grigori; M. Masha Sosonkina (2008): {\it A partitioning algorithm for block-diagonal matrices with overlap},  Parallel Computing,
{\bf 34}, 332--344.

\bibitem{karp1976} R. M. Karp (1976): {\it The probabilistic analysis of some combinatorial search algorithms},
in Algorithms and Complexity: New Directions and Recent Results,  1--20.


 \bibitem{KHO} S. Khot;  R. Saket (2008): {\it  Hardness of Minimizing and Learning DNF Expressions}, in
Proc. FOCS, pp. 231--240.


\bibitem{KRI} M. Krivelevich (1997): {\it Approximate set covering in uniform hypergraphs}. J. Algorithms
{\bf 25} , no. 1, pp. 118--143.

\bibitem{lan2007} G. Lan (2007): {\it An effective and simple heuristic for the set covering problem}, European Journal
of Operational Research, {\bf 176}, 1387-1403.


\bibitem{levi2008} A. Levin (2008): {\it  Approximating the unweighted k-set cover problem: greedy meets local
search}, SIAM J. Discrete Math., {\bf  231}, 25--264.

\bibitem{TAO2005} T. Li (2005): {\it  A general model for clustering binary data}, in Proceedings of the $11^th$ ACM SIGKDD int. conf. on Knowledge discovery in Data Mining (KDD '05). ACM, New York, NY, USA, 188--197.

\bibitem{LOV} L. Lovasz (1975): {\it On the ratio of the optimal integral and fractional covers}. Disc. Math.
{\bf 13}, pp. 383--390.

\bibitem{LUND} C. Lund; M. Yannakakis (1994): {\it On the hardness of approximating minimization problems},  J. ACM {\bf 31} , no. 5, pp. 960--981.

\bibitem{mezard1987} M. Mezard; G.  Parisi; M. A. Virasoro (1987): {\it Spin glass theory and beyond}, World Scientific, Singapore.

\bibitem{OKUN} M. Okun (2005): {\it On the approximation of the vertex cover problem in hypergraphs}. Discrete Optimization {\bf 2}, no. 1, pp. 101--111.

\bibitem{RAZ} R. Raz; M. Safra (2007): {\it A sub-constant error-probability low-degree test, and a subconstant error-probability PCP characterization of NP}. In Proc. STOC, pp.
475--484.

\bibitem{SAK} R. Saket; M. Sviridenko (2012): {\it New and Improved Bounds for the Minimum Set Cover Problem}, Approximation, Randomization, and Combinatorial Optimization.  Algorithms and Techniques
Lecture Notes in Computer Science, Volume 7408, pp 288--300.


\bibitem{SAM} A. Samorodnitsky; L. Trevisan (2000): {\it A PCP characterization of NP with optimal
amortized query complexity}, in Proc. STOC, pp. 191--199.


\bibitem{VE1984} C. Vercellis (1984):  {\it A Probabilistic Analysis of the Set Covering Problem}, Annals of Operations  Research  {\bf 1},  255--271.

\bibitem{VIJE2010} C.N. Vijeyamurthy; R. Panneerselvam (2010):
{\it Literature review of covering problem in operations management},
International Journal of Services, Economics and Management, {\bf 2}, 267--285.

\bibitem{ZHANG2005} Z.  Zhang; T. Li; C. Ding; X. Zhang (2007):  {\it Binary Matrix Factorization with Applications}, in Proceedings of the 2007 Seventh IEEE International Conference on Data Mining (ICDM '07).
IEEE Computer Society, Washington, DC, USA, 391--400.


}

\end{thebibliography}
\end{document}